\documentclass[journal]{IEEEtran}

\usepackage{array}
\usepackage{amsmath}
\usepackage{amsfonts}
\usepackage{amssymb}
\usepackage{mathrsfs}
\usepackage{latexsym}
\usepackage[dvips]{epsfig}
\usepackage[font=footnotesize]{subfig}
\usepackage{hyperref}
\usepackage{graphics}

\usepackage{enumerate}

\usepackage{amsthm}
\usepackage{arydshln}
\usepackage{cite}

\usepackage{listings}
\usepackage{color} %red, green, blue, yellow, cyan, magenta, black, white
\usepackage{verbdef}% http://ctan.org/pkg/verbdef
\definecolor{mygreen}{RGB}{28,172,0} % color values Red, Green, Blue
\definecolor{mylilas}{RGB}{170,55,241}
\usepackage{csquotes}
\usepackage{bibentry}
\usepackage{mathtools}
\usepackage{url}  %Required

\usepackage{amsmath,amssymb,amsfonts,mathrsfs,bm} 
\usepackage{caption}
\usepackage{cleveref}
\usepackage{array}
\usepackage{scalerel}
\usepackage{algorithm}
\usepackage{booktabs}
\usepackage{multirow}
\usepackage{rotating}
\usepackage{algpseudocode}
\usepackage{cite}
\usepackage{epsfig}
\usepackage{tabularx}
\usepackage{epstopdf}
\usepackage{circuitikz}
\usepackage{textcomp}
\usepackage{pgfplots}
\usepackage{lscape}
\pgfplotsset{compat=1.14}
\usepackage{array,booktabs}
\usepackage{siunitx}

\usepackage{lipsum}
\usepackage{mathtools}
\usetikzlibrary{arrows,shapes,positioning}
\usepackage{xspace}

\ifx\theorem\undefined
\newtheorem{theorem}{Theorem}
\newenvironment{theorem*}{\par\noindent{\bf Theorem\ }}{\hfill\\[2mm]}

\newenvironment{corollary*}{\par\noindent{\bf Corollary\ }}{\hfill\\[2mm]}
\newtheorem{definition}{Definition}

\fi

\newcommand{\Abf}{\boldsymbol{A}}

\newcommand{\bbf}{\boldsymbol{b}}

\newcommand{\cbf}{\boldsymbol{c}}

\newcommand{\Dbf}{\boldsymbol{D}}

\newcommand{\ebf}{\boldsymbol{e}}

\newcommand{\Fbf}{\boldsymbol{F}}

\newcommand{\Ibf}{\boldsymbol{I}}

\newcommand{\Kcal}{\mathcal{K}}

\newcommand{\Lbf}{\boldsymbol{L}}
\newcommand{\Lcal}{\mathcal{L}}

\newcommand{\Nbb}{\mathbb{N}}
\newcommand{\Ncal}{\mathcal{N}}

\newcommand{\Rbb}{\mathbb{R}}

\newcommand{\Ubf}{\boldsymbol{U}}

\newcommand{\Vbf}{\boldsymbol{V}}

\newcommand{\wbf}{\boldsymbol{w}}

\newcommand{\xbf}{\boldsymbol{x}}

\newcommand{\ybf}{\boldsymbol{y}}

\newcommand{\zbf}{\boldsymbol{z}}

\newcommand{\zerobf}{\boldsymbol{0}}

\newcommand{\eps}{\varepsilon}

\definecolor{Fcolor}{rgb}{0, 0.5, 0.25}
\newif\ifcomment

\crefrangelabelformat{equation}{(#3#1#4)\,--\,(#5#2#6)}
\crefmultiformat{equation}{(#2#1#3)}{\,--\,(#2#1#3)}{#2#1#3}{#2#1#3}
\crefname{equation}{}{}

\crefmultiformat{table}{(#2#1#3)}{\,--\,(#2#1#3)}{#2#1#3}{#2#1#3}
\crefname{table}{}{}%
\crefname{figure}{Figure}{Figures}
\crefname{algorithm}{Algorithm}{}
\crefname{table}{Table}{Tables}
\crefname{lemma}{Lemma}{Lemmas}
\crefname{theorem}{Theorem}{Theorems}
\crefname{section}{Section}{Sections}
\crefname{definition}{Definition}{Definitions}
\commentfalse

\makeatletter
\DeclareRobustCommand{\cev}[1]{%
	\mathpalette\do@cev{#1}%
}
\newcommand{\do@cev}[2]{%
	\fix@cev{#1}{+}%
	\reflectbox{$\m@th#1\vec{\reflectbox{$\fix@cev{#1}{-}\m@th#1#2\fix@cev{#1}{+}$}}$}%
	\fix@cev{#1}{-}%
}
\newcommand{\fix@cev}[2]{%
	\ifx#1\displaystyle
	\mkern#2 1mu
	\else
	\ifx#1\textstyle
	\mkern#2 3mu
	\else
	\ifx#1\scriptstyle
	\mkern#2 2mu
	\else
	\mkern#2 2mu
	\fi
	\fi
	\fi
}

\graphicspath{{_images/}}
\begin{document}
\title{A First-Order Numerical Algorithm\\ without Matrix Operations}
%
% not built to handle multiple paragraphs

\author{Muhammad Adil, Ramtin Madani, Sasan Tavakkol, and Ali Davoudi% <-this % stops a space
\thanks{Muhammad Adil is with Palo Alto Research Center (PARC), Ramtin Madani, and Ali Davoudi are with the University of Texas at Arlington. Sasan Tavakkol is with Google Research. Ramtin Madani (ramtin.madani@uta.edu) is the corresponding author.	
This work is funded, in part, by the Office of Naval Research under award N00014-18-1-2186, and approved for public release under DCN\# 43-9295-22.}}

\maketitle

\IEEEpeerreviewmaketitle

\begin{abstract}

This paper offers a matrix-free first-order numerical method to solve large-scale conic optimization problems. Solving systems of linear equations pose the most computationally challenging part in both first-order and second-order numerical algorithms. Existing direct and indirect methods are either computationally expensive or compromise on solution accuracy. Alternatively, we propose an easy-to-compute decomposition method to solve sparse linear systems that arise in conic optimization problems. Its iterations are tractable, highly parallelizable, with closed-form solutions. This algorithm can be easily implemented on distributed platforms, such as graphics processing units, with orders-of-magnitude time improvement. The performance of the proposed solver is demonstrated on large-scale conic optimization problems and is compared with the state-of-the-art first-order solvers. 
\end{abstract}


\section{Introduction}

Conic optimization is used in various areas such as operation research, machine learning, signal processing, and optimal control. Common solutions to conic optimization problems are based on the interior point method (IPM) \cite {Meh92, DCB13, TTT03, NW06, NY18, SW97}, that are suitable for small- to medium-sized problems. At each iteration, IPM solves a linear system of equations, mainly using Gauss-Jordan \cite{AM87}, Gaussian elimination \cite{SP13}, LU decomposition \cite{MK02}, Cholesky decomposition \cite{ST12}, QR decomposition, or Monte Carlo Methods \cite{FB07}. These direct approaches become prohibitively expensive by scale, and, rendering IPM-based methods impractical for large-scale problems.

Matrix-free interior-point methods, including indirect or iterative methods, are among the most popular methods for solving large-scale problem \cite{diamond2016matrix,JG09}. However, direct methods become prohibitive \cite{diamond2016matrix,TK04} at larger scale and Krylov subspace iterative methods, such as preconditioned conjugate gradient (PCG), become attractive alternatives \cite{TK04}. In conjugate gradient method \cite{HS52}, an iterative technique solves the Newton step instead of factorizing the Hessian matrix directly. \cite{FKS02} proposed a Lagrangian dual predictor-corrector algorithm, and applied the conjugate gradient method to solving linear systems. In \cite{KS07}, an iterative solver is applied to the modified barrier method for large-scale semidefinite programming optimization problems. The conjugate gradient method, with a simple preconditioner, reduces the computation time of solving semidefinite programming problems \cite{CY00}, and an inexact semismooth Newton conjugate gradient approach in \cite{ZST10} improves the solution accuracy.

The performance of iterative methods depends on the spectral properties of the linear system and the condition number of the matrix involved \cite{K95}. Existing preconditioners aim at reducing the condition number to the extend possible \cite{MVD10}. 
A matrix-free algorithm for equality-constraint nonlinear programming in \cite{CNW09} remedies ill-conditioned problems with rank-deficient Jacobian matrices.  A matrix-free solver PDCO \cite{S18} uses least squares minimal residual (LSMR) to solve linear systems. \cite{JG12} proposes a matrix-free IPM for quadratic programs, where the Karush–Kuhn–Tucker (KKT) system is regularized to bound the condition number and, then, a preconditioner is designed for the regularized system. 
This approach substantially decreases the computational cost in each iteration, albeit by trading off more iterations and lower accuracy. 

Collectively, IPMs are inherently computationally expensive and hard to scale for larger conic optimization problems. Alternatively, first-order methods scale gracefully with moderate accuracy \cite{OCPB16, ZFP+20, FB2018, SBGB2018}. The computational complexity of each iteration is significantly less than IPM. They are 
suitable when high accuracy is not needed \cite{BPC+11, SMP19}. The convergence of first-order methods in a limited number of iterations is an active area of research \cite{ASM21, NLR+15, GB17,EY15, BG18, BG16, FZB19, HL2017, OPY+2020, PC11}.  

First-order methods are required to solve a one time linear system \cite{SBGB2018}. Thanks to \textit{factorizing caching}, in practice, the linear system is solved only in the first iteration, factors are stored, and, then, reused in subsequent iteration \cite{SBGB2018, FB2018, SMP19}. The most common factorizing caching approaches are LDL and QR decompositions \cite{CONP13}. Direct methods or factorizing caching become impractical at a very large scale, and point matrix-free or indirect approaches become viable \cite{FB2018}. The first-order methods \cite{CONP13, OCPB16, FB2018, SBGB2018} apply conjugate gradient method to the resulting large-scale problems. First-order methods struggle with accuracy, and matrix-free indirect conjugate gradient approaches further hamper their ability. The computational complexity of direct methods and low accuracy solution of indirect methods to tackle large conic problems are the motivations behind the proposed matrix-free approach for solving very large problems with a modest accuracy. 
\subsection{Contributions}
We develop a matrix-free first-order numerical algorithm to solve very large-scale sparse conic optimization problems. The basic idea is to decompose the constraint matrix into sparse factors, such that iterative steps do not involve any matrix operations.  
These factors are easy-to-compute and require minimal memory storage. The matrix inversion lemma makes the operations matrix division free. We reformulate the standard conic optimization problem by introducing auxiliary variables and, then, apply the proposed matrix-free algorithm in conjunction with the well-known two-block alternating direction method of multipliers (ADMM) \cite{BPC+11, MKL15,MKL18}. Therefore, the computational burden of solving the linear system in each iteration is taken out of the iterative loop.  
The proposed algorithm admits parallel implementation, and is amenable to graphics processing units (GPUs). We demonstrate the performance gain and computational speedup of the proposed algorithm by conducting a range of experiments and compare the results with other first-order numerical solvers.  

\subsection{Paper Structure}

The rest of this paper is organized as follows. We introduce the cone programming and a brief description of operator splitting methods for such problems in Section \ref{sec:prelim}. We analyze the existing direct and indirect methods for solving the linear system, and provide a motivation for our numerical algorithm in Section \ref{sec: linear_sys}. Section \ref{sec: formulation} presents the proposed algorithm. Numerical experiments and comparison with competing solvers are presented in \ref{sec: experiments}. Section \ref{sec:conclusion} concludes the paper.

%\newpage 
%\ 
%\newpage
 \subsection{Notations}
Symbols $\Rbb$ and $\Nbb$ denote the set of real and natural numbers, respectively.  Matrices and vectors are represented by bold uppercase and bold lowercase letters, respectively. Notation $\lVert \cdot {\rVert}_2$ refers to $\ell_2$ norm of either matrix or vector depending on the context, and $\lvert \cdot\rvert$ represents the absolute value. The symbol $\!(\cdot)^{\!\top}\!$ represent the transpose operators. The notations $\boldsymbol{I}_n$ refer to the $n\times n$ identity  matrix. The symbol $\Kcal$ describes different types of cones. The superscript $\!(\cdot)^{\!\mathrm{opt}}\!$ refers to the optimal solution of optimization problem. $o$ denotes the number of nonzero entries in matrix. The symbol $\Lcal$ represent the augmented Lagrangian function. The notations ``$\mathcal{P}_1$", ``$\mathcal{P}_2$", represent the primal, whereas ``$\mathcal{D}$" denote the dual block of two-block ADMM. The symbols $\eps_{\rm prim} , \eps_{\rm dual}, \eps_{\mathrm{abs}}$, and $\eps_{\mathrm{rel}}$ are used for primal, dual, absolute, and relative tolerance, respectively. $\eps_{\rm gap}$ is the tolerance for difference between primal and dual objective values.

\section{Preliminaries}\label{sec:prelim}
We consider the class of convex optimization problems with linear objective function, linear constraints, and conic constraints, in the form of: 
\begin{subequations}\label{eq:prob_primal}
	\begin{align}
	& \underset{
		\begin{subarray}{c} \!\!\!\! \!\!\! \! \xbf\in\,\Rbb^{n}
		\end{subarray}
	}{\text{minimize~~~}}
	& &\hspace{-2cm} \cbf^{\top}\xbf  \label{eq:prob_obj}\\
	& \text{subject to~~~}
	& &\hspace{-2cm}  \Abf \xbf  = \bbf  \label{eq:prob_constraint}\\
	& & & \hspace{-2cm} \xbf \in \Kcal  \label{eq:prob_cone}
	\end{align}
\end{subequations}

where $\xbf\in\Rbb^n$ is the primal decision variable and $\cbf \in \Rbb^{n}$, $\Abf \in \Rbb^{m\times n}$, and $\bbf \in \Rbb^{m}$ are given. Additionally, $\Kcal \triangleq\Kcal_{n_1}\times\Kcal_{n_2}\times\cdots\times\Kcal_{n_k}\subseteq\Rbb^n$ is a non-empty, closed, convex cone, where each $\Kcal_{n_i}\subseteq\Rbb^{n_i}$ is a Lorentz cone of size $n_i$, i.e.,
\begin{align}
\Kcal_{n_i}\triangleq\big\{\wbf\in\Rbb^{n_i}\,|\,w_1\geq\big\|[w_2,\ldots,w_{n_i}]\big\|_2\big\},\nonumber
\end{align}
and $n_1+n_2+\ldots+n_k=n$.

In order to solve \ref{eq:prob_primal}, various first-order operator splitting methods have been proposed in the past decade. First-order methods are particularly interesting for the cases where iterative steps can be solved efficiently through explicit formula, and a large number of iterations can be executed in a short amount of time. Among most popular methods are the Douglas-Rachford Splitting (DRS) technique and ADMM. %To this end, we provide a description of DRS and ADMM methods to solve problem \ref{eq:prob_primal}. 

\subsubsection{Douglas-Rachford Splitting}
DRS was originally proposed in \cite{DR56} to find numerical solutions of differential equations for heat conduction problems, and it has been widely used to solve separable convex optimization problems. Rather than operating on the whole problem directly, DRS works on a splitting scheme to address each component of the problem separately.  
In order to implement the DR splitting method, one casts the problem \eqref{eq:prob_primal} in the form of
\begin{align}
&{\text{minimize}}
&&\hspace{-2cm} f(\xbf) + g(\xbf)  \label{eq:split_opt}
\end{align}
where $f,g:\Rbb^n\to\Rbb\cup\{\infty\}$ %and $g:\Rbb^n\to\Rbb\cup\{\infty\}$ 
are the indicator functions:
%are proper closed and convex functions. 
\begin{align}
f(\xbf)\! \triangleq\!  \left\{\begin{matrix} 0 & \!\!\mathrm{if }\  \xbf \in \Kcal \\
\infty & \!\!\mathrm{ \ otherwise} \end{matrix}\right.  
\quad\mathrm{and}\quad
g(\xbf) \!\triangleq\!  \left\{\begin{matrix} \cbf^{\top}\xbf & \!\!\mathrm{ if }\ \Abf\xbf=\bbf \\ \infty & \!\!\mathrm{otherwise} \end{matrix}\right.  \nonumber 
\end{align}
leading to the following iterative steps:
\begin{subequations}
	\begin{align}
	&\xbf\gets\mathrm{prox}_f(\zbf)\\
	&\zbf\gets\zbf+\mathrm{prox}_{\mu^{-1}g}(2\xbf-\zbf)-\xbf
	\end{align}
\end{subequations}
where $\mu$ is a fixed tuning parameter, and for every $h:\Rbb^n\to\Rbb\cup\{\infty\}$ and $\wbf_0\in\Rbb^n$, the operator $\mathrm{prox}_h(\wbf_0)$ returns the unique solution to the following optimization problem:
%\begin{subequations}
\begin{align}
& \underset{
\begin{subarray}{c} \!\!\!\! \!\!\! \! \wbf \in \Rbb^n
\end{subarray}
}{\text{minimize~~~}}
&&\hspace{-2cm} h(\wbf)+\frac{1}{2}\|\wbf-\wbf_0\|_2^2.
\end{align}
%\end{subequations}
Each iteration of DRS requires the evaluation of the proximal operators $\mathrm{prox}_f(\cdot)$ and $\mathrm{prox}_{\mu^{-1}g}(\cdot)$. While the evaluation of $\mathrm{prox}_f(\cdot)$ is parallelizable and enjoys a closed-form solution, the evaluation of $\mathrm{prox}_{\mu^{-1}g}(\wbf_0)$ is the main bottleneck which requires solving the following system of linear equations:
\begin{align}
\begin{bmatrix}
\Ibf_n&\Abf^{\top}\\
\Abf&\zerobf
\end{bmatrix}
\begin{bmatrix}
\wbf\\
\boldsymbol{\nu}
\end{bmatrix}=
\begin{bmatrix}
\wbf_0-\mu^{-1}\cbf\\
\bbf
\end{bmatrix},\label{lin}
\end{align}
where $\wbf=\mathrm{prox}_{\mu^{-1}g}(\wbf_0)$.
\vspace{2mm}

\subsubsection{Alternating Direction Method of Multipliers} ADMM is one of the most commonly used first-order methods for solving large-scale optimization problems. %ADMM is a method of choice for large-scale optimization problems owing to its simple implementation, computational cheapness, and distributive nature.  
ADMM can be analyzed as a special case of DRS, as the former applied to the primal problem is equivalent to the latter applied to the dual problem. A standard way of solving problem \eqref{eq:prob_primal} via ADMM is through the following reformulation
\begin{subequations}
	\begin{align}
	& \underset{
		\begin{subarray}{c} \!\!\!\! \!\!\! \! \xbf_1,\xbf_2 \in \Rbb^n
		\end{subarray}
	}{\text{minimize~~~}}
	&&\hspace{-2cm} f(\xbf_1) + g(\xbf_2) \\
	& \text{subject to~~~}
	&&\hspace{-2cm} \xbf_1=\xbf_2 \label{xx}
	\end{align}
\end{subequations}
leading to the steps
\begin{subequations}
	\begin{align} \label{eq:admm_p1}
	&\xbf_1 \gets \mathrm{prox}_{f}(\xbf_2-\mu^{-1}\boldsymbol{\lambda_1})\\ \label{eq:admm_p2}
	&\xbf_2 \gets \mathrm{prox}_{\mu^{-1}g}(\xbf_1+\mu^{-1}\boldsymbol{\lambda_1})\\ 
	&\boldsymbol{\lambda}_{12}\gets\boldsymbol{\lambda}_{12}+\mu(\xbf_1-\xbf_2).
	\end{align}
\end{subequations}
where $\xbf_1, \xbf_2$ are primal variables, $\boldsymbol{\lambda}_{12}$ is the dual variable associated with \eqref{xx}, and $\mu$ is a fixed tuning parameter. Each iteration requires solving the system of linear equations \eqref{lin}.

\section{Solving the Linear System} \label{sec: linear_sys}
While first-order methods are considered highly efficient for solving large-scale problems with modest accuracy, evaluating the projection operator $\mathrm{prox}_{\mu^{-1}g}$ at each iteration can become computationally prohibitive. 
Depending on the structure and size of the constraint matrix $\Abf$, there are different ways to solve the linear system \eqref{lin}. The common approaches for solving \eqref{lin} in several iterations are classified as direct and indirect methods \cite{FB2018,SBGB2018}. 

\subsubsection{Direct Methods}

The linear system \eqref{lin} can be solved by first factoring the matrix in \eqref{lin} and, then, performing forward or backward substitutions. The common practice is to compute the factors at the beginning and, then, reuse them in subsequent iterations. This process %of computing the factors at the first iteration and then reusing in rest of the iterations 
is known as \textit{factorization caching}, %. The technique of factorization caching is very effective to attain a large speedup in subsequent iterations, after the first iteration and is frequently 
which is widely used by  first-order solvers \cite{FB2018, SBGB2018, OCPB16, CONP13}.
When $\Abf$ is dense, the common practice is to solve the linear system \eqref{lin} by performing either Cholesky or QR decompositions and, then, cache the factors for subsequent iterations. The complexity in both approaches is $O(mn^2)$ in first iteration and $O(mn)$ in subsequent iterations. When $\Abf$ is sparse and not too large, the direct sparse $\Lbf\Dbf\Lbf^{\top}$ decomposition of the matrix in \eqref{lin} can be used 
with $\Lbf$ being a lower triangular matrix and $\Dbf$ being a diagonal matrix \cite{Van95}. %The solution of the linear system can be made division free by storing $\Dbf^{-1}$ in the first step. 
In this case, the computational complexity %for sparse matrices 
%is $O(o_{\Lbf})$ and   
depends upon the number of nonzero entries of the factor $\Lbf$.    

\subsubsection{Indirect Methods}

Instead of solving \eqref{lin} by factorization, \textit{indirect methods} apply an iterative procedure such as conjugate gradient (CG) \cite{HS52,GV96, NS06} or LSQR \cite{PS82, HDWC13}.   
%Indirect methods are required to evaluate the approximate projections in each step,  but  in practice,  this method takes very few iterations to produce an adequate approximation. A more practical approach is to terminate the method when the linear system is solved  to a reasonable accuracy. Another advantage of indirect method is that the conjugate gradient technique can be warm-started by using the solution obtained in previous iteration as an initial point in each subsequent step. This heuristic approach results  inaccurate projections at the beginning but start producing more accurate ones as the first-order algorithms begins to converge.  
In contrast to direct methods, indirect methods offer lower complexity and  more freedom for adaptive parameter selection since there is no factorization. In each step of conjugate gradient, a vector is multiplied by $\Abf$ and $\Abf^{\top}$, thus, the computational complexity of each step is significantly lower. % $O(o_{\Abf})$. 

\section{UV Decomposition} \label{sec: formulation}

In contrast to existing direct and indirect methods, herein, we propose an efficient matrix-free decomposition to handle large-scale sparse optimization problems of the form \eqref{eq:prob_primal}. 
\begin{definition}
%Let $o$ represent the number of non-zero elements in $\Abf\in\Rbb^{m\times n}$ and 
Assume that $\Abf=\Ubf\Vbf^\top$, where $\Ubf \in \Rbb^{m\times o}$ and $\Vbf \in \Rbb^{n \times o}$. The pair of matrices $(\Ubf,\Vbf)$ are regarded as a UV decomposition of the matrix $\Abf$, if $\Ubf\Ubf^\top$ and $\Vbf\Vbf^\top$ are both diagonal.
\end{definition}
 
Such decomposition can be readily constructed as follows. 
Let $\{\Abf_{i_k,j_k}\}^o_{k=1}$ represent the non-zero elements of $\Abf$ in an arbitrary order. The pair
\begin{subequations}
	\begin{align}
	&\Ubf\triangleq\sum_{k=1}^{o}{\Abf_{i_k,j_k}\ebf_{i_k}\dot{\ebf}_k^{\top}}\label{UU}\\
	&\Vbf\triangleq\sum_{k=1}^{o}{\dot{\ebf}_k^{\top}\ddot{\ebf}_{j_k}}\label{VV}
	\end{align}
\end{subequations}
is a possible candidate for a UV decomposition of $\Abf$,
where $\{\ebf_k\}^m_{k=1}$, $\{\dot{\ebf}_k\}^o_{k=1}$, and $\{\ddot{\ebf}_k\}^n_{k=1}$ represent the standard basis for $\Rbb^m$,  $\Rbb^o$, and $\Rbb^n$, respectively. 
The factor $\Ubf$ in \eqref{UU} contains the nonzero elements of the matrix $\Abf$, whereas $\Vbf$ in \eqref{VV} encodes the locations. Next, we further elucidate this decomposition by providing a simple illustrative example.
\vspace{2mm}

\noindent\textbf{Example 1:} %We illustrate the decomposition $\Abf=\Ubf\Vbf$, where $\Abf$ is a sparse matrix with $o$ nonzero entries, the corresponding $\Ubf$ and $\Vbf$ are 
Consider the following sparse matrix with $8$ non-zero elements:
\begin{subequations}
\begin{align} 
\Abf &:= \begin{bmatrix} 1  &  0 &   4 &   6 &    8 \\
						 0  &  0 &   5 &   0 &    0 \\
						 2  &  3 &   0 &   7 &    0 \\
\end{bmatrix}.
\end{align}
The pair
\begin{align}
\Ubf &:= \begin{bmatrix} 1  &  0 &   0 &   4 &    0 &    6 &    0 &    8 \\
						 0  &  0 &   0 &   0 &    5 &    0 &    0 &    0 \\
						 0  &  2 &   3 &   0 &    0 &    0 &    7 &    0\\
\end{bmatrix},\\
\Vbf &:= \begin{bmatrix} 1  &  1 &   0 &   0 &    0 &   0 &   0 &    0      \\
						 0  &  0 &   1 &   0 &    0 &   0 &   0 &    0     \\
						 0  &  0 &   0 &   1 &    1 &   0 &   0 &    0     \\
						 0  &  0 &   0 &   0 &    0 &   1 &   1 &    0     \\
						 0  &  0 &   0 &   0 &    0 &   0 &   0 &    1     \\
\end{bmatrix}
\end{align}
\end{subequations}
satisfies $\Abf=\Ubf\Vbf^{\top}$. Additionally, each column of $\Ubf$ and each column of $\Vbf$ has exactly one non-zero element. Hence, $(\Ubf,\Vbf)$ is a UV decomposition of the matrix $\Abf$.

Next, using the notion of UV decomposition, we seek to develop a matrix-free numerical algorithm for solving problem \ref{eq:prob_primal}. To this end, we reformulate the problem as follows:
\begin{subequations}\label{eq:prob_reform}
	\begin{align}
	& \underset{
		\begin{subarray}{c} \!\!\!\! \!\!\! \! \xbf\in\,\Rbb^{n}\\
			\!\!\! \!\!\! \! \ybf\in\,\Rbb^{o}\\
			\!\!\!\!\! \!\! \! \zbf\in\,\Rbb^{n}
		\end{subarray}
	}{\text{minimize~~~}}
	& & \hspace{-2cm}\cbf^{\top}\xbf \\
	& \text{subject to~~~}\label{eq:prob_U}
	& & \hspace{-2cm} \boldsymbol{\lambda}:\;\Ubf \ybf  = \bbf\\ \label{eq:prob_V}
	& & & \hspace{-2cm}\boldsymbol{\gamma}:\;\ybf=\Vbf^{\top}\xbf\\   \label{eq:prob_z}
	& & & \hspace{-2cm}\boldsymbol{\delta}:\;\zbf=\xbf\\ 
	& & & \hspace{-2cm}\zbf \in \Kcal\hspace{2cm}
	\end{align}
\end{subequations}
where $\ybf\in\Rbb^o$ and $\zbf\in\Rbb^n$ are auxiliary variables, and $(\Ubf,\Vbf)$ is a UV decomposition of the matrix $\Abf$. 
It is straightforward to verify that the two formulations are equivalent.
To solve problem \eqref{eq:prob_reform} via ADMM, we form the augmented Lagrangian function:
%\begin{equation}
\begin{align}
\!\!\!\Lcal(\xbf,\ybf, \zbf, \boldsymbol{\lambda},\boldsymbol{\gamma},\boldsymbol{\delta})\triangleq
\cbf^{\top}\xbf
& + \boldsymbol{\lambda}^{\top}(\Ubf\ybf - b) +\dfrac{\mu}{2}\big\|\Ubf\ybf-\bbf\big\|_2^2 \nonumber\\
&+ \boldsymbol{\gamma}^{\top}(\ybf\!-\!\Vbf^{\!\top}\xbf) + \dfrac{\mu}{2}\big\|\ybf\!-\!\Vbf^{\!\top}\xbf\big\|_2^2\nonumber\\
&+ \boldsymbol{\delta}^{\top}(\zbf-\xbf) + \dfrac{\mu}{2}\big\|\zbf-\xbf\big\|_2^2\label{lll}
\end{align}	
%\end{equation}
where $\mu>0$ is a fixed parameter, and $\boldsymbol{\lambda} \in \Rbb^{m}, \boldsymbol{\gamma} \in \Rbb^{o}$, and $\boldsymbol{\delta}\in \Rbb^{n}$ are the Lagrange multipliers associated with constraints \ref{eq:prob_U},\ref{eq:prob_V} and \ref{eq:prob_z}, respectively.

We perform two block ADMM and regroup the primal and dual variables as follows
\begin{align*}
\hspace{-2cm}\text{(Block 1)} \hspace{1cm} \mathcal{P}_1 &= \{\xbf\}\\
\hspace{-2cm}\text{(Block 2)} \hspace{1cm}\mathcal{P}_2 &= \{\ybf,\zbf\}\\
\hspace{-2cm}\text{(Dual)} \hspace{1.4cm}\mathcal{D}   &= \{\boldsymbol{\lambda}, \boldsymbol{\gamma}, \boldsymbol{\delta}\} 
 \end{align*} 

Given the above partitioning of variables, each iteration of two block ADMM for problem \ref{eq:prob_reform} involves the following steps:

%\begin{enumerate}
%	\item 
\subsubsection{Block 1} This step consist of minimizing the Lagrangian function \eqref{lll} with respect to the variable $\xbf$, and freezing the other variables at their previous values, i.e., 
	\begin{align*}
	\xbf^{k+1}=
	\underset{
		\begin{subarray}{c} \!\!\!\! \!\!\! \! \xbf \in \Rbb^n
		\end{subarray}
	}{\text{\rm {arg\ min~~~}}}
	\Lcal(\xbf,\ybf^{k}, \zbf^{k}, \boldsymbol{\lambda}^{k},\boldsymbol{\gamma}^{k},\boldsymbol{\delta}^{k})
	\end{align*}
which enjoys the following closed-form solution
	\begin{align}\label{eq:vvt}
		\!\!\!\!\xbf^{k+1}\!=\!
	(\Ibf\!+\!\Vbf\Vbf^{\top})^{-1}\Big[
	\Vbf\big(\ybf^{k}\!+\!\frac{\boldsymbol{\gamma}^{k}}{\mu}\big)+
	\zbf^{k}+\frac{\boldsymbol{\delta}^{k}}{\mu}-\frac{\cbf}{\mu}\Big].\!\!\!
	\end{align}
	Observe that $(\Ibf+\Vbf\Vbf^{\top})$ is a diagonal positive definite matrix and, once the factor $\Vbf$ is known, computing the solution of \ref{eq:vvt} can be made matrix-free by storing the diagonal elements of $(\Ibf+\Vbf\Vbf^{\top})^{-1}$.  
%	\item

\subsubsection{Block 2}
The next block consists of minimizing the Lagrangian function \eqref{lll} with respect to the variables $\ybf$ and $\zbf$, while $\xbf$ is fixed:
\begin{align*}
	(\ybf^{k+1},\zbf^{k+1})=
	\underset{
		\begin{subarray}{c} \!\!\!\! \!\!\! \! \ybf \in \Rbb^o\\
			\!\!\!\! \!\! \! \zbf \in \Rbb^n
		\end{subarray}
	}{\text{\rm {arg\ min~~~}}}
	\Lcal(\xbf^{k+1},\ybf, \zbf, \boldsymbol{\lambda}^{k},\boldsymbol{\gamma}^{k},\boldsymbol{\delta}^{k})
\end{align*}
which involves two parallel steps.

\textit{Minimization with respect to $\ybf$}: leads to the following closed-form solution: 
\begin{align}\!\!\!\ybf^{k+1}\!=\!
(\Ibf\!+\!\Ubf^{\top}\Ubf)^{-1}\Big[
\Ubf^{\top}\big(\bbf-\frac{\boldsymbol{\lambda}^{k}}{\mu}\big)\!+\! 
\Vbf^{\top}\xbf^{k+1}-\frac{\boldsymbol{\gamma}^{k}}{\mu}\Big].
\end{align}
Now, according to the matrix inversion lemma, we have
\begin{align}
(\Ibf+\Ubf^{\top}\Ubf)^{-1}=(\Ibf- \Ubf^{\top}(\Ibf + \Ubf\Ubf^{\top})^{-1}\Ubf)
\end{align}
and since $(\Ibf + \Ubf\Ubf^{\top})$ is a diagonal positive-definite matrix, this step can be made matrix free as well, by storing the diagonal elements of $(\Ibf + \Ubf\Ubf^{\top})^{-1}$.

%		\item 
\textit{Minimization with respect to $\zbf$}: involves projection onto the associated Lorentz cones which is parallelizable and enjoys a closed-form solution as well: 
\begin{align}
\zbf^{k+1}=\mathrm{proj}_{\Kcal}\Big\{\xbf^{k+1}-\frac{\boldsymbol{\delta}^{k}}{\mu}\Big\}.
\end{align}
Observe that for each $i\in\{1,\ldots,k\}$,
\begin{align} \label{eq:proj}
\mathrm{proj}_{\Kcal_i}(\wbf)  \triangleq  \left\{\begin{matrix} \hspace{-2.5cm}\mathbf{0}^{n_i} & \!\!\alpha\leq  -w_1 \\
\hspace{-2.7cm}\wbf &\!\! \alpha\leq  +w_1 \\
\frac{\wbf}{2}\!+\![\frac{\alpha}{2}, \frac{w_1w_2}{2\alpha},\ldots, \frac{w_1w_{n_i}}{2\alpha}]  & \!\!\mathrm{ \ otherwise} \end{matrix}\right. 
\end{align}
where $\alpha=\|[w_2,\ldots,w_{n_i}]\|_2$.

%	\end{enumerate} 
%	\item 
	
\subsubsection{Dual variables update} This steps involves the update of dual variable as follows: 
\begin{subequations}
\begin{align}
\boldsymbol{\lambda}^{k+1} &= \boldsymbol{\lambda}^{k}+\mu(\Ubf \ybf^{k+1}-\bbf)\\
\boldsymbol{\gamma}^{k+1} &= \boldsymbol{\gamma}^{k}+\mu(\ybf^{k+1}-\Vbf^{\top}\xbf^{k+1})\\
\boldsymbol{\delta}^{k+1} &= \boldsymbol{\delta}^{k}+\mu(\zbf^{k+1}-\xbf^{k+1})
\end{align}
\end{subequations}

%\end{enumerate}

\begin{algorithm}[t]
	\caption{~}
	\label{alg:1}
	\vspace{0.2mm}
	\begin{algorithmic}[1]
		\Require{ $(\Abf,\bbf,\cbf,\Kcal)$, fixed $\mu > 0$, and initial points $\xbf, \zbf, \boldsymbol{\delta}\in\Rbb^n,  \ybf,\boldsymbol{\gamma} \in \Rbb^o$, and $\boldsymbol{\lambda} \in \Rbb^m$}
		\State Construct a  UV decomposition of the matrix $\Abf$
		\State $\Fbf_{\Ubf} := \big(\Ibf + \Ubf\Ubf^{\top}\big)^{-1}$
		\State $\Fbf_{\Vbf}:= \big(\Ibf + \Vbf\Vbf^{\top}\big)^{-1}$
		
		\Repeat
		%\State $k \gets k+1$
		\State $ \xbf \gets  \Fbf_{\Vbf}\Big(
		\Vbf\big(\ybf+\frac{\boldsymbol{\gamma}}{\mu}\big)+
		\zbf+\frac{\boldsymbol{\delta}}{\mu}-\frac{\cbf}{\mu}\Big)$
		\State $ \ybf \gets  (\Ibf - \Ubf^{\top}\Fbf_{\Ubf}\Ubf)\Big(
		\Ubf^{\top}\big(\bbf-\frac{\boldsymbol{\lambda}}{\mu}\big)+ \Vbf^{\top}\xbf-\frac{\boldsymbol{\gamma}}{\mu}\Big)$
		\State $\zbf \gets \mathrm{proj}_{\Kcal}\Big\{\xbf-\frac{\boldsymbol{\delta}}{\mu}\Big\}$
		
		\State			$\boldsymbol{\lambda} \gets \boldsymbol{\lambda}+\mu(\Ubf \ybf-\bbf)$
		\State $\boldsymbol{\gamma} \gets \boldsymbol{\gamma}+\mu(\ybf-\Vbf^{\top}\xbf)$
		\State $\boldsymbol{\delta} \gets \boldsymbol{\delta}+\mu(\zbf-\xbf)$
		
		\Until {stopping criteria is met.$\phantom{\Big|}$}
		\Ensure \!
		$\xbf^{\mathrm{opt}}, \boldsymbol{\lambda}^{\mathrm{opt}}$

	\end{algorithmic}\label{al:alg_1}
\end{algorithm}
%In contrary to existing practices that focuses on either standard matrix decompositions or apply an approximate solution by conjugate gradient method, %we propose a matrix free algorithm which does not rely on standard matrix factorization. Instead, 
%we propose to store the nonzero entries in sparse matrices in such a way that we do not need to apply the standard matrix decomposition. The standard fast ADMM iterations can be used to the reformulated problem. 
The above steps are summarized in Algorithm \ref{al:alg_1}. 

\begin{itemize}
	\item {\bf Steps 1:} The equality constraint matrix $\Abf$ is decomposed into $\Abf = \Ubf\Vbf^{\top}$ such that $\Ubf\Ubf^{\top}$ and $\Vbf\Vbf^{\top}$ are both diagonal. 
	\item {\bf Steps 2 and 3:}  To make iterative steps matrix-free, we compute the multiplication factors in these steps. These factors are easy-to-compute since $(\Ibf + \Ubf\Ubf^{\top})$ and $(\Ibf + \Vbf\Vbf^{\top}
	)$  are positive definite diagonal matrices. 
	\item {\bf Steps 5 to 10:} These steps are the ADMM primal and dual variables updates, as well as the projection onto specific cones performed in Step 7.   
\end{itemize}

\section{Numerical Experiments} \label{sec: experiments}
We highlight the computational strength and scalability of the proposed algorithm by testing it on a variety of randomly generated large-scale and sparse linear programming (LP) and second-order cone programming (SOCP) problems. We compare our methods with the first-order solvers POGS (Proximal Graph Solver) \cite{FB2018}, OSQP (Operator Splitting Solver for Quadratic Programs \cite{SBGB2018}, and SCS (Splitting Conic Solver) \cite{OCPB16}.
For each experiment, we stop immediately after a better solution than the competing solver is obtained in terms of residual norms, constraint violations, or objective gap. The performance gain of the proposed algorithm is consistent across all problem instances and different hardware architectures.    

The proposed algorithm and competing solvers are implemented in MATLAB R2020a, and all the experiments are carried out  on a Linux-based DGX station with 20 2.2 GHz cores, Intel Xeon E5-2698 v4 CPU, with NVIDIA Tesla V100-DGXS-32GB (128 GB total) GPU processor and 256 GB of RAM. The parallel nature of the proposed algorithm helps us to take advantage of multi-core CPU. Note that our implementation in MATLAB utilizes only a single GPU and does not benefit from multiple GPUs. Moreover, no experiment is bounded by RAM or GPU memory of the DGX station. We have used the MATLAB interface of POGS, OSQP v0.6.0 and SCS v2.1.2.

\textit{Problem instances:} All data is generated in such a way that the linear and second order cone programming problems are feasible and bounded. The number of nonzero entries of $\Abf$ are in the range of $10^4$ to $10^7$. 
\begin{itemize}
	\item   We generate $\Abf \in \Rbb^{m\times n}$ to be a sparse random matrix with $0.1\%$, $0.5\%$,  and $1.0\%$ nonzero elements, which are drawn i.i.d. (independently and identically distributed) from $\Ncal(0,1)$.
	\item $\bbf: =\Abf\times \mathrm{proj}_{\Kcal}(\dot{\xbf})$, where the elements of $\dot{\xbf}\in\Rbb^{n}$ have i.i.d standard normal distribution.
	\item The elements of $\cbf\in\Rbb^n$ have i.i.d standard normal distribution.
	\item For linear programming instances, we have $\Kcal = \Rbb_+^n$. For SOCP instance, we have $\Kcal=(\Kcal_4)^{\frac{n}{4}}$, where $\Kcal_4$ is the standard Lorentz cone of size $4$.
\end{itemize}

\subsection{Linear Programming}
We consider randomly-generated linear programming problems, and compare the performance of algorithm \ref{al:alg_1} with POGS and OSQP on a variety of sparse problems. Experiments are continued until the run time of the competing solver reaches a maximum time of 1200 seconds. The maximum time is chosen in such a way that the experiments provide sufficient information to compare the computational time for all solvers.

\begin{figure*}[t] 
	
	\subfloat[\label{fig:osqp_0p1}] {\includegraphics[width = 0.33\textwidth]{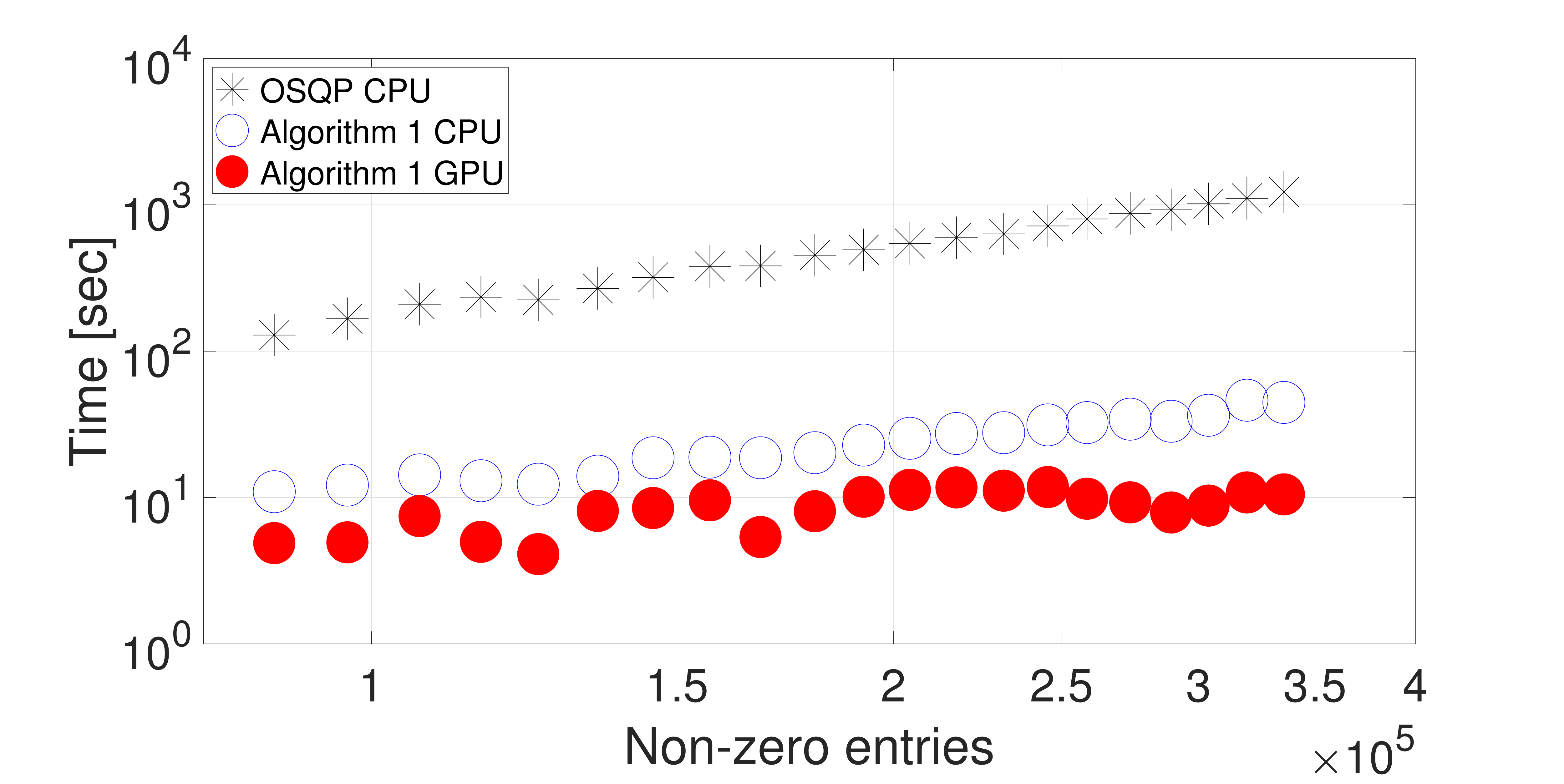}}
	\subfloat[\label{fig:osqp_0p5}] {\includegraphics[width = 0.33\textwidth]{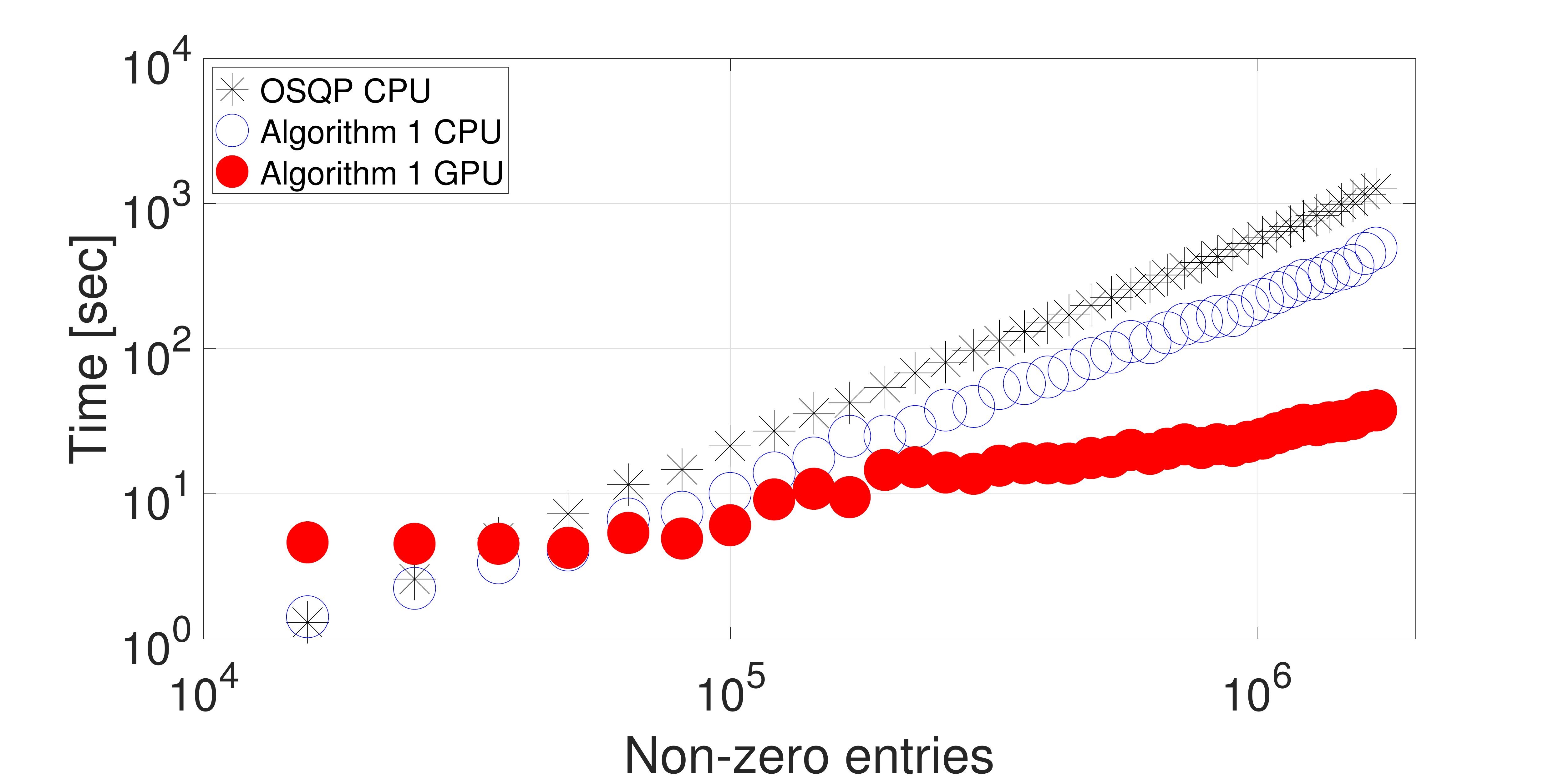}}
	\subfloat[\label{fig:osqp_1p0}] {\includegraphics[width = 0.33\textwidth]{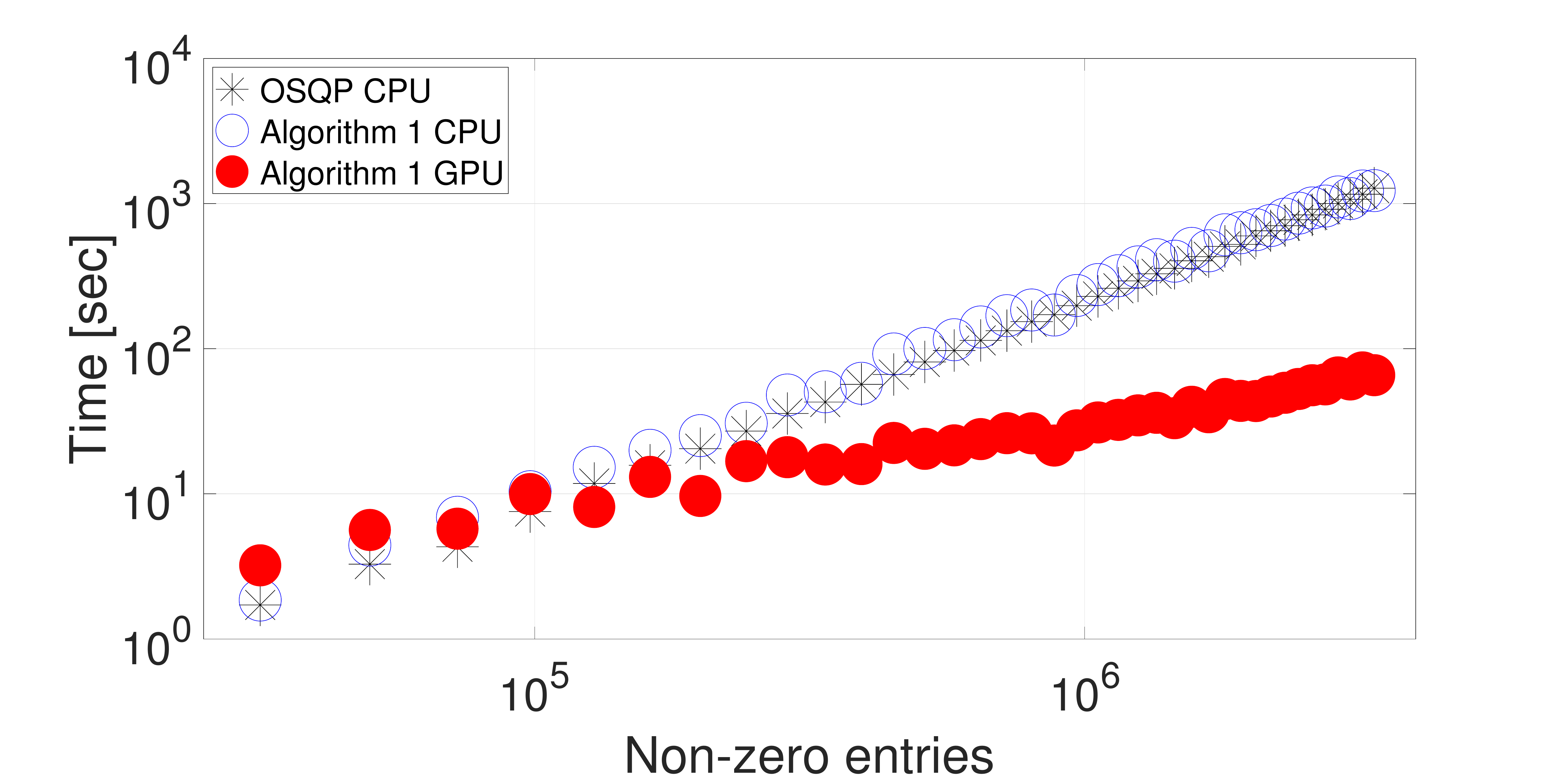}}	
	\caption{The performance of Algorithm \ref{al:alg_1}  for linear programming in comparison with  OSQP for (a) $0.1\%$ (b) $0.5\%$ and (c) $1.0\%$ nonzero entries in matrix $\Abf$.} 
	
	\label{fig:lp_computation_time_op1}
\end{figure*}

\begin{figure*}[t] 
	
	\subfloat[\label{fig:pogs_0p1}] {\includegraphics[width = 0.33\textwidth]{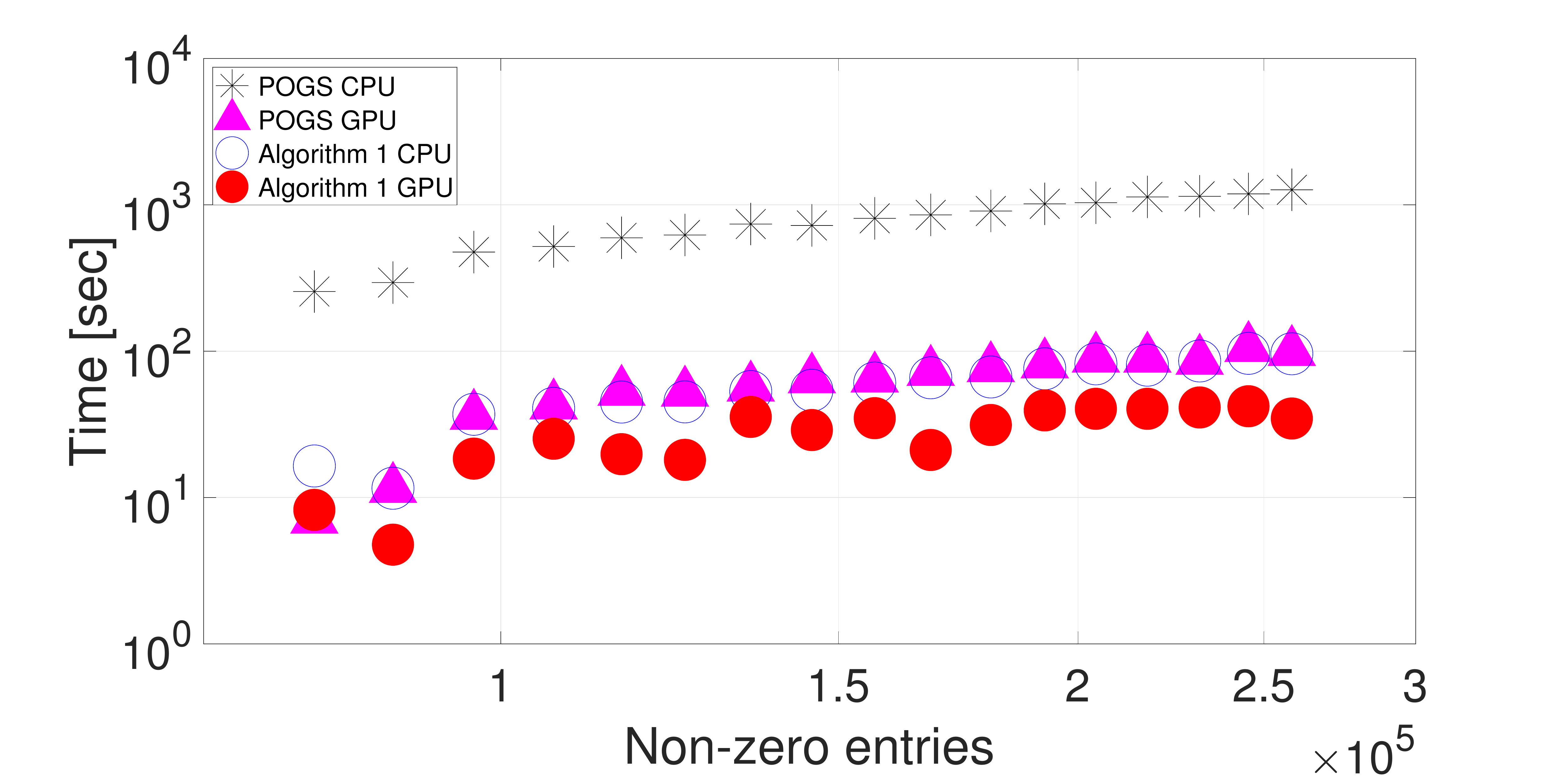}}
	\subfloat[\label{fig:pogs_0p5}] {\includegraphics[width = 0.33\textwidth]{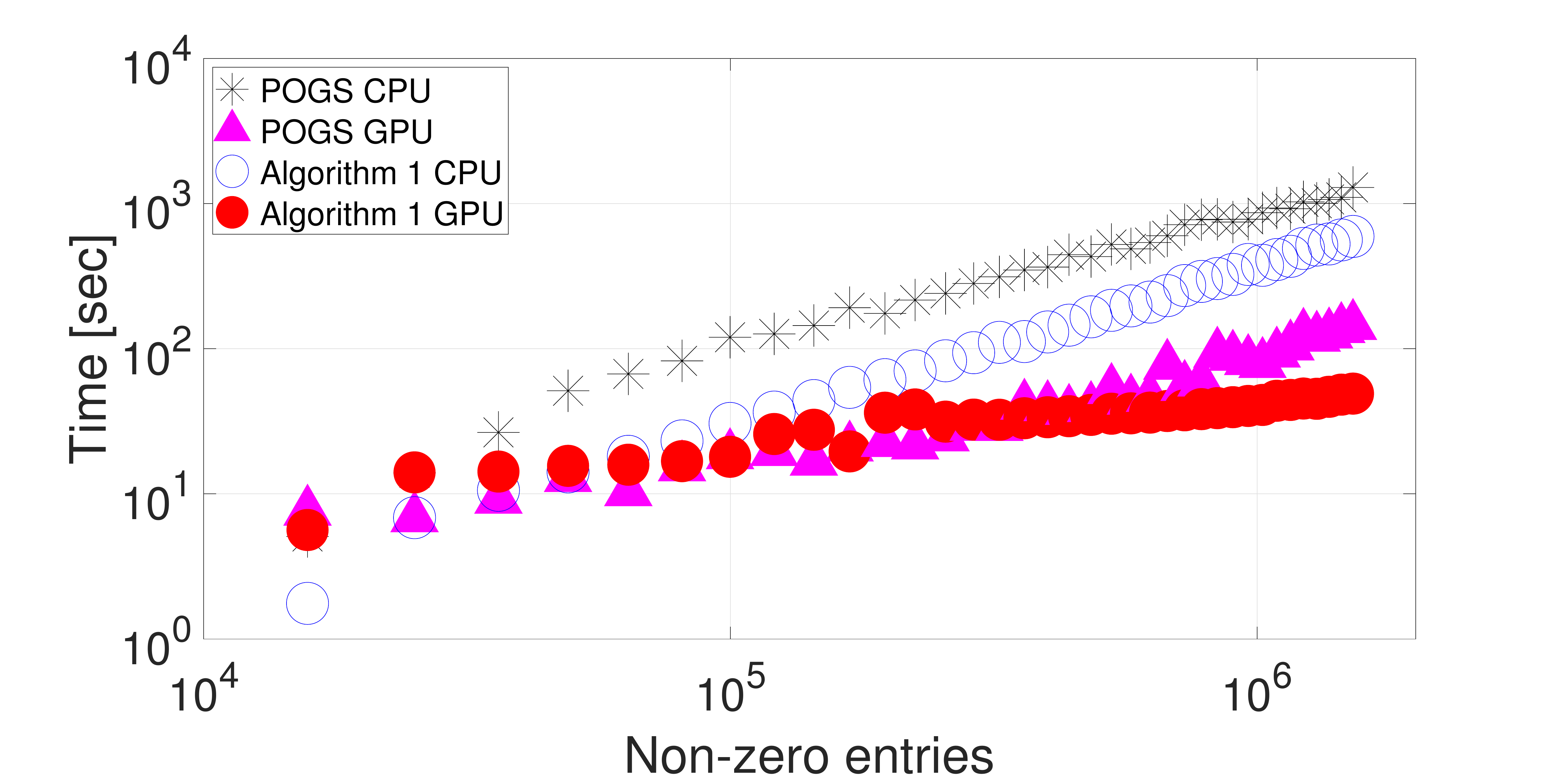}}
	\subfloat[\label{fig:pogs_1p0}] {\includegraphics[width = 0.33\textwidth]{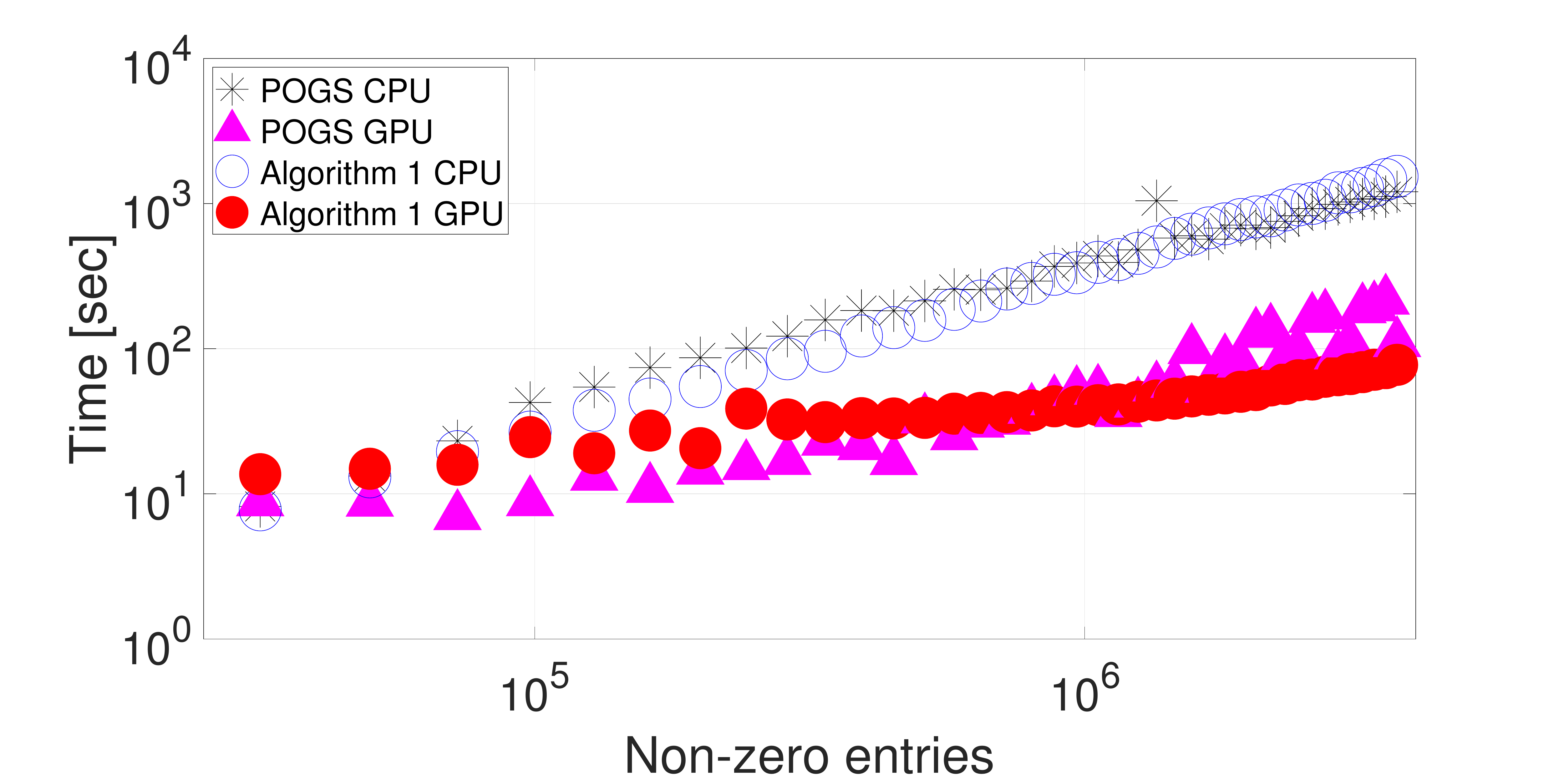}}
	\caption{The performance of Algorithm \ref{al:alg_1}  for linear programming in comparison with  POGS for (a) $0.1\%$  (b) $0.5\%$ and (c) $1.0\%$ nonzero entries in matrix $\Abf$. }

	\label{fig:lp_pogs}
\end{figure*}

\subsubsection{Comparison with OSQP}

OSQP is a first-order general purpose open-source solver based on the alternating direction method of multipliers. In Fig. \ref{fig:lp_computation_time_op1}, we compare our method with OSQP by running several sparse linear programming instances. For each instance, we use OSQP solver in its default settings and, then, provide these parameters as input to our algorithm to achieve the same tolerance values.
\newline
\textit{Termination criteria:} We stop when both primal and dual residuals are smaller than some predefined tolerance limits $\varepsilon_{\rm prim} > 0$ and $\varepsilon_{\rm dual} > 0$, i.e.,
\begin{subequations}\label{stopcri_osqp}
\begin{align}
&\|\Abf \xbf-\bbf\|_{\infty} < \varepsilon_{\rm prim}\\
& \|\Abf^{\top} \boldsymbol{\lambda} + \cbf\|_{\infty} < \varepsilon_{\rm dual}\\
&\eps_{\rm prim} = \eps_{\rm abs} + \eps_{\rm rel} \max\lbrace \|\Abf \xbf\|_{\infty}, \| \bbf \|_{\infty} \rbrace \\
&\eps_{\rm dual} = \eps_{\rm abs} + \eps_{\rm rel} \max\lbrace  \| \Abf^\top \boldsymbol{\lambda} \|_{\infty}, \| \cbf \|_{\infty} \rbrace
\end{align}
\end{subequations}
where  $\eps_{\rm abs} = 10^{-4}$ and $\eps_{\rm rel} = 10^{-3}$ are  default absolute and relative tolerance values, respectively. We also satisfy the following condition to terminate our algorithm 
\begin{align*}
| \cbf^{\top}\xbf  + \bbf^{\top} \boldsymbol{\lambda}|<| \cbf^{\top}\xbf^{\rm OSQP} + \bbf^{\top}\ybf^{\rm OSQP} |
\end{align*}

 Figure \ref{fig:osqp_0p1} demonstrates that the CPU implementation of the proposed algorithm is, at least, ten times faster than OSQP, whereas the GPU implementation shows two orders-of-magnitude improvement. Similar results are depicted in Figures \ref{fig:osqp_0p5} and \ref{fig:pogs_1p0}, where we have used $0.5\%$ and $1.0\%$ nonzero values for the constraint matrix $\Abf$. The CPU and GPU implementations of the proposed algorithm can easily achieve ten times and hundred times improvements, respectively, while satisfying similar, or even stricter, stopping criteria as compared with OSQP.

\subsubsection{Comparison with POGS}

 We compare the performance of the proposed algorithm with POGS on both CPU and GPU system architectures in figure \ref{fig:lp_pogs}. POGS is an open-source implementation of graph projection splitting method that targets multi-core and GPU-based systems for solving convex optimization problems. Herein, POGS would return inaccurate solution in its default parameter settings. Hence, we made a slight change in its default parameters by setting $\eps_{\mathrm{abs}} = 10^{-5}$ and $\eps_{\mathrm{rel}} = 10^{-4}$ to obtain reasonable accuracy. The same tolerance parameters, along with primal and dual solutions returned by POGS, are used as input parameters for the proposed algorithm to meet the same stopping criteria.
\newline
\textit{Termination criteria:} The stopping criteria of Algorithm \ref{alg:1} is when we exceed both primal and dual feasibility of the solution produced by the competing solver, i.e.,, when the following two criteria are met: 
\begin{subequations}\label{stopcri}
	\begin{align}
	&\|\Abf \xbf-\bbf\|_2 < \|\Abf \xbf^{\rm POGS}-\bbf\|_2\\
	&| \cbf^{\top}\xbf + \bbf^{\top} \boldsymbol{\lambda}|<|  \cbf^{\top}\xbf^{\rm POGS} + \bbf^{\top}\boldsymbol{\lambda}^{\rm POGS}|
	\end{align}
\end{subequations}
where $\xbf^{\text{POGS}}$ and $\boldsymbol{\lambda}^{\text{POGS}}$ are primal and dual solutions produced by the competing solver POGS under slightly-modified default settings.

Figure \ref{fig:pogs_0p1}, \ref{fig:pogs_0p5}, and \ref{fig:pogs_1p0}, show notable time improvements offered by the proposed algorithm in both cases. Our CPU and GPU implementation of Algorithm \ref{alg:1} shows an order of magnitude time improvements compared to POGS, and this improvement is further pronounced for larger problems. 

\subsection{Second-order Cone Programming}

We compare the performance of Algorithm \ref{al:alg_1} in with splitting conic solver (SCS) on a variety of sparse conic problems.  

\subsubsection{Comparison with SCS} Primarily written in C, SCS a first-order numerical optimization solver for large-scale cone programs. This solver returns both primal and dual solutions along with infeasibility certificate when applies.
\newline
\textit{Termination criteria:} We use the following default stopping criteria of the competing solver 
\begin{subequations}
\begin{align}
&\|\Abf\xbf - \bbf\|_2 \leq \eps_{\rm prim}(1+\|\bbf\|_2),\\ &\|\Abf^{\top}\boldsymbol{\lambda} + \cbf\|_2 \leq \eps_{\rm dual}(1+\|\cbf\|_2),\\
&\cbf^{\top}\xbf + \bbf^{\top}\boldsymbol{\lambda} \leq \eps_{\rm gap}(1+|\cbf^{\top}\xbf| + |\bbf^{\top}\boldsymbol{\lambda}|)
\end{align}
\end{subequations}
where $\eps_{\rm prim}=\eps_{\rm dual}= \eps_{\rm gap} = 10^{-3}$.

\begin{figure*}[t]  
	
	\subfloat[\label{fig:scs_0p1}] {\includegraphics[width = 0.33\textwidth]{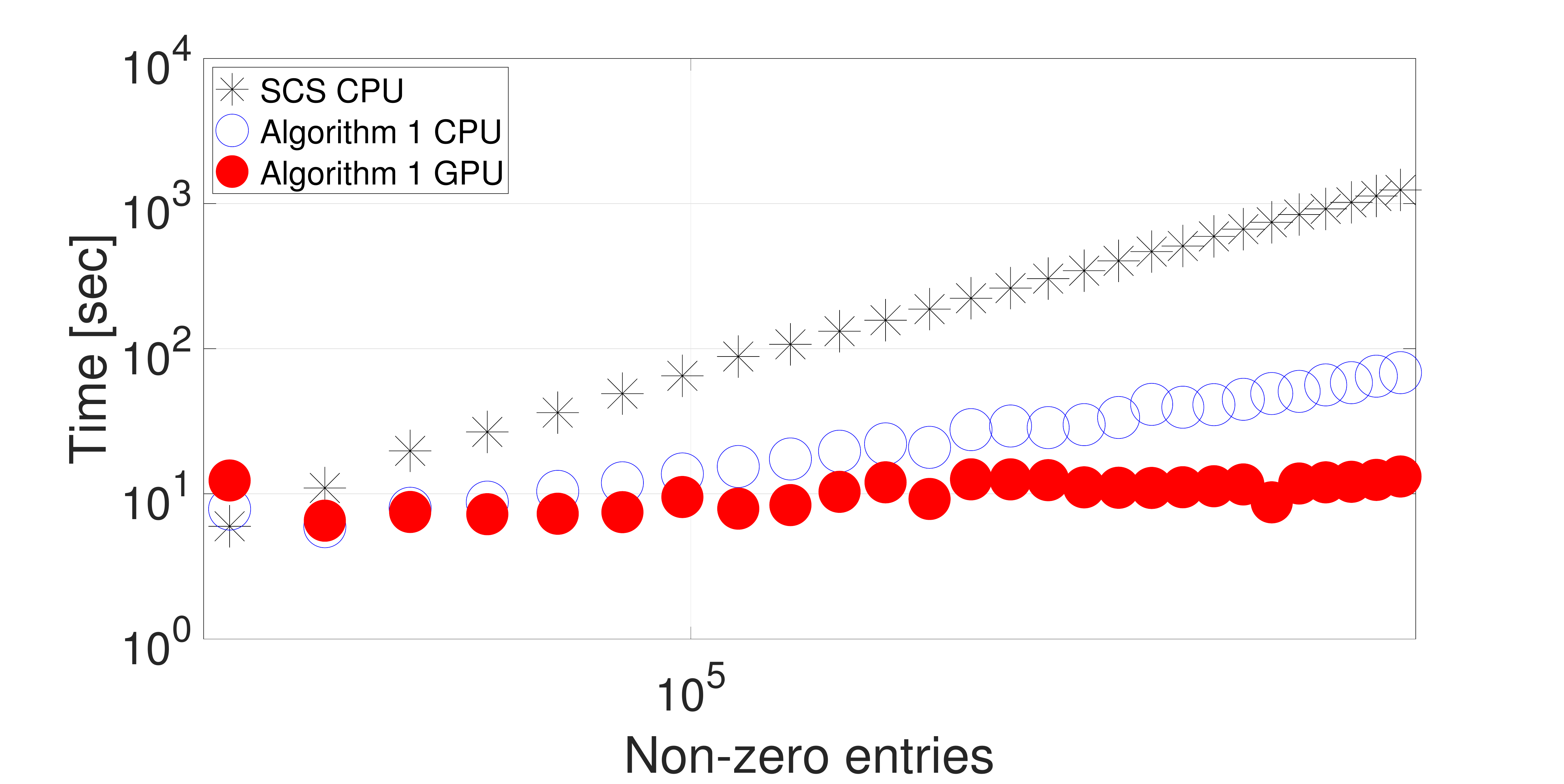}}
	\subfloat[\label{fig:scs_0p5}] {\includegraphics[width = 0.33\textwidth]{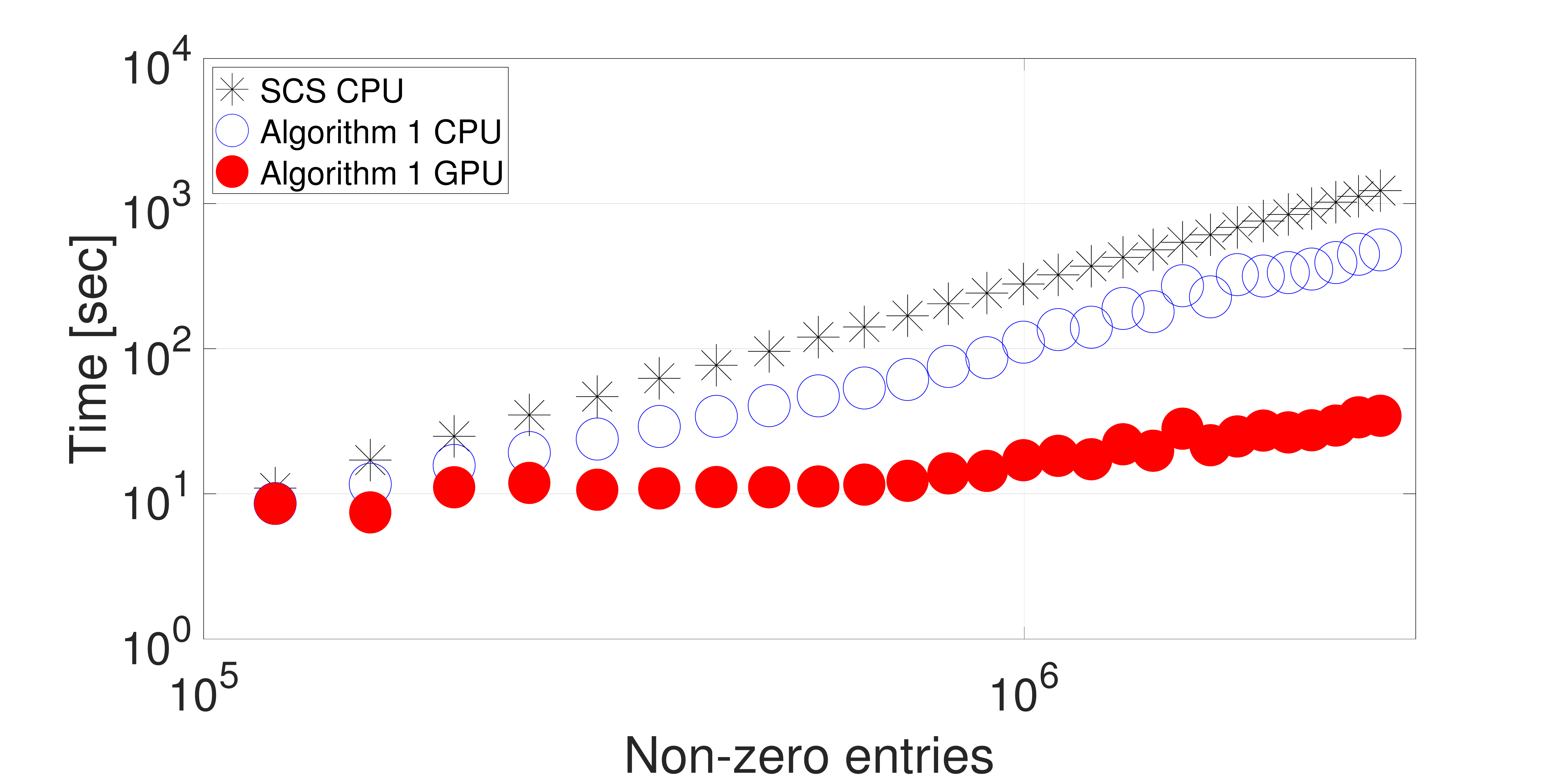}}
	\subfloat[\label{fig:scs_1p0}] {\includegraphics[width = 0.33\textwidth]{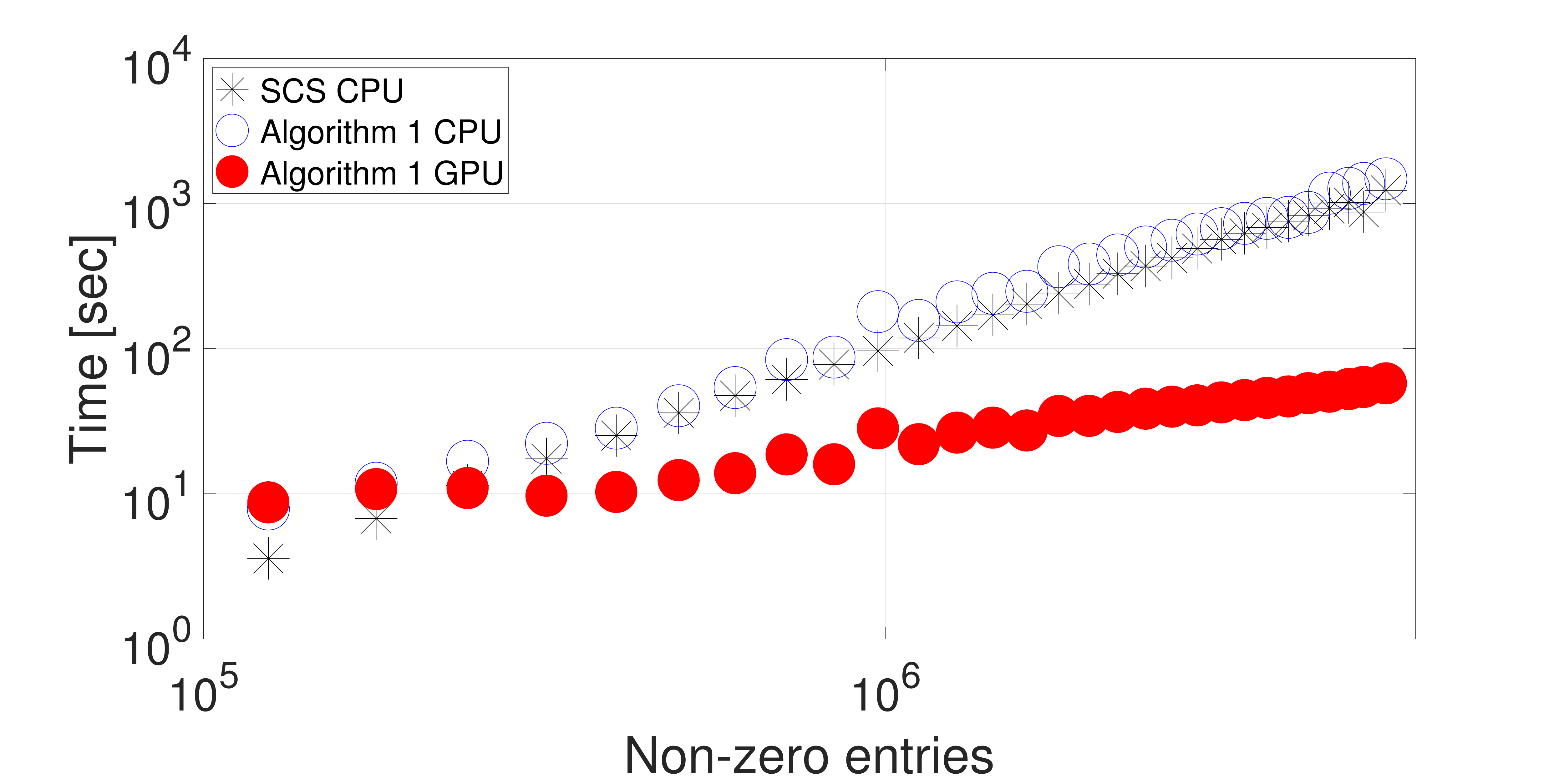}}
	
\caption{The performance of Algorithm \ref{al:alg_1} for second order cone programming  in comparison with  SCS direct  for (a) $0.1\%$ (b) $0.5\%$ and (c) $1.0\%$ nonzero entries in matrix $\Abf$. }
	
	\label{fig:socp}
\end{figure*}

We have solved multiple SOCP problems for different sparse density and satisfy the termination criteria of SCS solver and compare the computational time in Figure \ref{fig:socp}. 

\section{Conclusion} \label{sec:conclusion}

We propose a computationally-efficient matrix-free first-order method for solving large-scale sparse conic optimization problems.  The computational burden is managed by decomposing the equality constraint into sparse factors that are easy-to-compute.  
A highly parallelizable and computationally-cheap numerical algorithm is developed based on the alternating direction method of multipliers. Each iteration of the proposed algorithm has a closed-form solution with simple arithmetic operations. The parallel nature of the proposed algorithm allows graphics processing unit (GPU) implementation and speeds up the computational gains by an order-of-magnitude. The proposed algorithm is applied to several linear and conic optimization problems. The numerical experiments show that the proposed matrix-free algorithm significantly achieves approximately an order-of-magnitude time improvement in comparison with POGS, OSQP, or SCS solvers.

%\vspace{0.5\textheight}
%
%{\ }

\bibliographystyle{IEEEtran}
\bibliography{IEEEabrv,egbib_revised}

%\appendix

\end{document}